\input amstex
\documentstyle{amsppt}

\define\I{\text{\rm{Im}}\,}



\topmatter
\title
On the  P\'olya-Wiman properties of Differential Operators
\endtitle
\author
Min-Hee Kim and Young-One Kim
\endauthor

\date
2015. 6. 1.
\enddate

\abstract Let $\phi(x)=\sum \alpha_n x^n$ be a formal power series
with real coefficients, and let $D$ denote differentiation. It is
shown that ``for every real polynomial $f$ there is a positive
integer $m_0$ such that $\phi(D)^mf$ has only real zeros whenever
$m\geq m_0$'' if and only if ``$\alpha_0=0$ or $2\alpha_0\alpha_2 -
\alpha_1^2 <0$'', and that if $\phi$ does not represent a
Laguerre-P\'olya function, then there is a Laguerre-P\'olya function
$f$ of genus $0$ such that  for every positive integer $m$,
$\phi(D)^mf$ represents a real entire function having infnitely many
nonreal zeros.
\endabstract

\address
Department of Mathematical Sciences, Seoul National University,
Seoul 151-747, Korea
\endaddress

\email alsgml01\@snu.ac.kr
\endemail

\address
Department of Mathematical Sciences and Research Institute of
Mathematics, Seoul National University, Seoul 151-747, Korea
\endaddress

\email kimyone\@snu.ac.kr
\endemail



\thanks
This work was supported by the SNU Mathematical Sciences Division
for Creative Human Resources Development.
\endthanks

\endtopmatter

\document

\head
1. Introduction
\endhead

A real entire function is an entire function which takes real values
on the real axis. If $f$ is a real entire function, we denote the
number of nonreal zeros (counting multiplicities) of $f$ by
$Z_C(f)$. (If $f$ is identically equal to $0$, we set $Z_C(f)=0$.) A
real entire function $f$ is said to be of {\it genus} $1^*$ if it
can be expressed in the form
$$
f(x)=e^{-\gamma x^2} g(x),
$$
where $\gamma\geq 0$ and $g$ is a real entire function of genus at
most $1$. (For the definition of the order and genus of an entire
function, see \cite{2, 14}.) If $f$ is a real entire function of
genus $1^*$ and $Z_C(f)=0$, then $f$ is called a {\it
Laguerre-P\'olya function} and we write $f\in\Cal{ LP}$. We denote
by $\Cal{ LP}^*$ the class of real entire functions $f$ of genus
$1^*$ such that $Z_C(f)<\infty$. It is well known that $f\in\Cal{
LP}$ if and only if there is a sequence $\langle f_n \rangle$ of
real polynomials such that $Z_C(f_n)=0$ for all $n$ and $f_n\to f$
uniformly on compact sets in the complex plane. (See Chapter 8 of
\cite{14} and \cite{15, 17, 21}.) From this and an elementary
argument based on Rolle's theorem, it follows that the classes
$\Cal{LP}$ and $\Cal{LP}^*$ are closed under differentiation, and
that $Z_C(f)\geq Z_C(f')$ for all $f\in\Cal{LP}^*$. The {\it
P\'olya-Wiman theorem} states that for every $f\in \Cal{ LP}^*$
there is a positive integer $m_0$ such that $f^{(m)} \in \Cal{ LP}$
for all $m\geq m_0$ \cite{7, 8, 10, 12, 20}. On the other hand, it
follows from recent results of  W. Bergweiler, A. Eremenko and J.
Langley that if $f$ is a real entire function, $Z_C(f)<\infty$ and
$f\notin\Cal{ LP}^*$, then $Z_C(f^{(m)})\to \infty$ as $m\to\infty$
\cite{1, 13}.

Let $\phi$ be a formal power series given by
$$
\phi(x) = \sum_{n=0}^{\infty} \alpha_n x^n.
$$
For convenience we express the $n$-th coefficient $\alpha_n$ of
$\phi$ as $\phi^{(n)}(0)/n!$ even when the radius of convergence is
equal to $0$. If $f$ is an entire function and the series
$$
\sum_{n=0}^{\infty} \alpha_n f^{(n)}
$$
converges uniformly on compact sets in the complex plane, so that it
represents an entire function, we write $f\in \text{\rm dom\,}
\phi(D)$ and denote the entire function by $\phi(D) f$. For $m\geq
2$ we denote by $\text{\rm dom\,} \phi(D)^m$ the class of entire
functions $f$ such that $f, \phi(D)f, \dots, \phi(D)^{m-1}f \in
\text{\rm dom\,} \phi(D)$. It is obvious that if $f$ is a
polynomial, then $f\in\text{\rm dom\,} \phi(D)^m$ for all $m$. For
more general restrictions on the growth of  $\phi$ and $f$ under
which $f\in\text{\rm dom\,} \phi(D)^m$ for all $m$, see \cite{4, 6}.

The following version of the  P\'olya-Wiman theorem for the operator
$\phi(D)$ was established by  T. Craven and G. Csordas \cite{6,
Theorem 2.4}.

\proclaim{Theorem A} Suppose that $\phi$ is a formal power series
with real coefficients, $\phi'(0)=0$ and $\phi''(0)\phi(0)<0$. Then
for every real polynomial $f$ there is a positive integer $m_0$ such
that all the zeros of $\phi(D)^m f$ are real and simple whenever
$m\geq m_0$.
\endproclaim

\remark{Remark} The assumption implies that $\phi(0)\ne 0$. On the
other hand, if $\phi(0)=0$ and $f$ is a real polynomial, then it is
trivial to see that $Z_C(\phi(D)^m f)\to 0$ as $m\to\infty$. (Recall
that we have set $Z_C( f)=0$ if $f$ is identically equal to $0$.)
\endremark

We also have the following version, which is a consequence of the
results in Section 3 of \cite{6}.

\proclaim{Theorem B} Suppose that $\phi\in \Cal{LP}$ ($\phi$
represents a Laguerre-P\'olya function), $f\in\Cal{LP}^*$, and that
 $f$ is of order less than  $2$. Then $f\in\text{\rm dom\,} \phi(D)^m$,
$\phi(D)^m f\in\Cal{LP}^*$ and $Z_C(\phi(D)^m f) \geq
Z_C(\phi(D)^{m+ 1}f)$ for all $m$. Furthermore, if $\phi$ is not of
the form $\phi(x)=ce^{\gamma x}$ with $c\ne 0$, then $Z_C(\phi(D)^m
f) \to 0$ as $m\to\infty$.
\endproclaim

\remark{Remarks} (1) If $\phi(x)=ce^{\gamma x}$, then
$$
\phi(D)f(x)=\sum_{n=0}^{\infty}\frac{c\gamma^n}{n!}f^{(n)}(x)=cf(x+\gamma)
$$
for every entire function $f$. Hence $Z_C(\phi(D)^mf)=Z_C(f)$ for
all $m$ whenever $c, \gamma \in\Bbb R$, $c\ne 0$ and $f$ is a real
entire function. We also remark that $\phi(x)=ce^{\gamma x}$ with
$c\ne 0$ if and only if $\phi(0)\ne 0$ and
$\phi^{(n)}(0)\phi(0)^{n-1}-\phi'(0)^n=0$ for all $n$.

(2) From \cite{6, Lemma 3.2}, \cite{11, Theorem 2.3} and the
arguments given in \cite{4}, it follows that the restriction ``$f$
is of order less than $2$'' can be weakened to ``$\phi$ or $f$ is of
genus at most $1$''. See also \cite {6, Theorem 3.3}.

(3) In the case where $\phi$ is of genus $2$, that is, $\phi$ is of
the form $\phi(x)=e^{-\gamma x^2}\psi(x)$, where $\gamma>0$ and
$\psi\in\Cal{LP}$ is of genus at most $1$, we have the following
stronger result: If $f$ is a real entire function of genus at most
$1$, and if the imaginary parts of the zeros of $f$ are uniformly
bounded, then $f\in\text{\rm dom\,} \phi(D)^m$ and $Z_C(\phi(D)^m f)
\geq Z_C(\phi(D)^{m+ 1}f)$ for all $m$, and $Z_C(\phi(D)^m f) \to 0$
as $m\to\infty$, even when $f$ has infinitely many nonreal zeros.
See \cite{3, Theorems 9a, 13 and 14},  \cite{4},  \cite{6, Lemma
3.2} and \cite{11, Theorem 2.3}.
\endremark

In this paper, we complement Theorems A and B above. Let $\phi$ be a
formal power series with real coefficients and $f$ be a real entire
function. If $f\in\text{\rm dom\,} \phi(D)^m$ for all $m$ and
$Z_C(\phi(D)^m f)\to 0$ as $m\to\infty$, then we will say that
$\phi$ (or the corresponding operator $\phi(D)$) has  the {\it
P\'olya-Wiman property} with respect to $f$. For instance, if $f$ is
a real entire function and $Z_C(f)<\infty$, then the operator $D$
($=d/dx$) has the P\'olya-Wiman property with respect to $f$ if and
only if $f\in\Cal{LP}^*$. Theorem A gives a sufficient condition for
$\phi$ to have the P\'olya-Wiman property with respect to arbitrary
real polynomials. The following two theorems imply that this is the
case if and only if $\phi(0)=0$ or $\phi''(0)\phi(0) - \phi'(0)^2
<0$.

\proclaim{Theorem 1.1} Suppose that $\phi$ is a formal power series
with real coefficients, $\phi(0)\ne 0$ and $\phi''(0)\phi(0) -
\phi'(0)^2 <0$. Then for every real polynomial $f$ there is a
positive integer $m_0$ such that all the zeros of $\phi(D)^m f$ are
real and simple whenever $m\geq m_0$.
\endproclaim

\proclaim{Theorem 1.2} Suppose that $\phi$ is a formal power series
with real coefficients, $\phi(0)\ne 0$,  $\phi''(0)\phi(0) -
\phi'(0)^2 \geq 0$, $\phi$ is not of the form $\phi(x)=ce^{\gamma
x}$ with $c\ne 0$,   $f$ is a real polynomial, and that
$$
\deg f \geq \min \{n\geq 2 :
\phi^{(n)}(0)\phi(0)^{n-1}-\phi'(0)^n\ne0\}.
$$
Then there  is a positive integer $m_0$ such that $Z_C(\phi(D)^m f)>
0$ for all $m\geq m_0$.
\endproclaim

If $\phi\in\Cal{LP}$ is not of the form $\phi(x)=ce^{\gamma x}$ with
$c\ne 0$, then it is easy to see that $\phi(0)=0$ or
$\phi''(0)\phi(0) - \phi'(0)^2 <0$ (for a proof, see \cite{5, 9});
hence Theorem B as well as Theorem 1.1 implies that $\phi$ has the
P\'olya-Wiman property with respect to arbitrary real polynomials.
On the other hand, there are plenty of formal power series $\phi$
with real coefficients which satisfy $\phi(0)=0$ or
$\phi''(0)\phi(0) - \phi'(0)^2 <0$, but do not represent
Laguerre-P\'olya functions. The following theorem, which is a strong
version of the converse of Theorem B, implies that if $\phi$ is one
of such formal power series, then $\phi$ does not have the
P\'olya-Wiman property with respect to some (transcendental)
Laguerre-P\'olya function of genus $0$, although it has the property
with respect to arbitrary real polynomials.

\proclaim{Theorem 1.3} Suppose that $\phi$ is a formal power series
with real coefficients and $\phi$ does not represent a
Laguerre-P\'olya function. Then there is a Laguerre-P\'olya function
$f$ of genus $0$ such that $f\in\text{\rm dom\,} \phi(D)^m$  and
$Z_C(\phi(D)^m f)=\infty$ for all positive integers $m$.
\endproclaim

As we shall see in the next section, Theorems 1.1 and 1.2 are almost
immediate consequences of Theorems 2.1 and 2.2 below, which are
proved in the same section by refining the arguments of Craven and
Csordas given in Section 2 of \cite{6}. Theorem 1.3 is a consequence
of P\'olya's characterization of the class $\Cal{LP}$ given in
\cite{18, 21} and a  diagonal argument. It is proved in Section 3.
Finally, in Section 4, we conclude the paper with some consequences
of Theorems 2.1 and 2.2 on  the asymptotic behavior of the
distribution of zeros of $\phi(D)^m f$ as $m\to\infty$, in the
general case where $\phi$ is a formal power series with complex
coefficients and $f$ is an arbitrary complex polynomial.

\head 2. Proof of Theorems 1.1 and 1.2
\endhead

For notational clarity, we denote the monic monomial of degree $d$
by $M^d$, that is, $M^d(x) = x^d$. With this notation, we have
$$
\left(\exp\left(\beta D^p\right)M^d\right)(x) = \sum_{k=0}^{\lfloor
d/p \rfloor} \frac{d! \beta^k}{k! (d-pk)!}x^{d-pk} \qquad
(\beta\in\Bbb C;\ d,p=1,2,\dots).
$$
The following two theorems will be proved after the proof of
Theorems 1.1 and 1.2.

\proclaim{Theorem 2.1} Suppose that $\phi$ is a formal power series
with complex coefficients, $\phi(0)=1$, $\phi$ is not of the form
$\phi(x)=e^{\gamma x}$,
$$
p=\min\{n : n\geq 2 \text{ \rm and }  \phi^{(n)}(0)\ne \phi'(0)^n\},
$$
$\alpha=\phi'(0)$ and
$\beta=\left(\phi^{(p)}(0)-\phi'(0)^p\right)/p!$. Suppose also that
$f$ is a monic complex polynomial of degree $d$, and $f_1, f_2,
\dots$ are given by
$$
f_m(x)=m^{-d/p} \left(\phi(D)^mf\right)(m^{1/p}x - m\alpha). \tag
2.1
$$
Then $f_m \to \exp\left(\beta D^p\right)M^d$ uniformly on compact
sets in the complex plane.
\endproclaim

\proclaim{Theorem 2.2} Suppose that $d$ and $p$ are positive
integers, $p\geq 2$, $q=\lfloor d/p \rfloor$ and $r=d-pq$. \roster
\item If $q=0$ ($d<p$), then $\exp\left(- D^p\right)M^d=M^d$.
\item If $q\geq 1$, then $\exp\left(- D^p\right)M^d$ has exactly
$q$ distinct positive zeros; and if we denote them by $\rho_1,
\dots, \rho_q$, then
$$
\left(\exp\left(- D^p\right)M^d\right)(x) = x^r\prod_{j=1}^{q}
\prod_{k=0}^{p-1}\left(x-e^{2k\pi i/p} \rho_j \right).
$$
\endroster
\endproclaim

\remark{Remark} The $d$-th Hermite polynomial $H_d$ is given by
$$
H_d(x) = \sum_{k=0}^{\lfloor d/2 \rfloor} \frac{(-1)^k d!}{k!
(d-2k)!}(2x)^{d-2k}.
$$
Thus we have $\left(\exp\left(- D^2\right)M^d\right)(x)=H_d(x/2)$
for all $d$, and Theorem 2.2 implies the well known fact that all
the zeros of the Hermite polynomials are real and simple.
\endremark

\proclaim{Corollary} If $\beta>0$, then all the zeros of
$\exp\left(-\beta D^2\right)M^d$ are real and simple, and
$\exp\left(\beta D^2\right)M^d$ has exactly $2 \lfloor d/2 \rfloor$
distinct purely imaginary zeros; and if $\beta\ne 0$ and $3\leq
p\leq d$, then $\exp\left(\beta D^p\right)M^d$ has nonreal zeros.
\endproclaim

This corollary is an immediate consequence of Theorem 2.2 and the
following relations which are trivially proved: If $\beta>0$ and
$\rho^p = -1$, then
$$
\left(\exp\left(-\beta D^p\right)M^d\right)(x) =
\beta^{d/p}\left(\exp\left(-D^p\right)M^d\right)\left(\frac{x}{\beta^{1/p}}\right)
$$
and
$$
\left(\exp\left(\beta D^p\right)M^d\right)(x)
=\left(\rho\beta^{1/p}\right)^d\left(\exp\left(-D^p\right)M^d\right)\left(\frac{x}{\rho\beta^{1/p}}\right).
$$

\demo{Proof of Theorem 1.1} Let $f$ be a (nonconstant) real
polynomial. Since multiplication by a nonzero constant does not
change the zeros of a polynomial, we may assume that $f$ is monic
and  $\phi(0)=1$. Let $d=\deg f$, $\alpha=\phi'(0)$,
$\beta=\left(\phi''(0)-\phi'(0)^2\right)/2$, and $f_1, f_2, \dots$
be given by
$$
f_m(x)=m^{-d/2}\left(\phi(D)^mf\right)(m^{1/2}x - m\alpha).  \tag
2.2
$$
Then $\beta<0$, and Theorem 2.1 implies that $f_m \to
\exp\left(\beta D^2\right)M^d$ uniformly on compact sets in the
complex plane. We have $\deg f_m=d=\deg(\exp\left(\beta
D^2\right)M^d)$ for all $m$; and since $\beta<0$, the corollary to
Theorem 2.2 implies that all the zeros of $\exp\left(\beta
D^2\right)M^d$ are real and simple. Hence the intermediate value
theorem implies that there is a positive integer $m_0$ such that all
the zeros of $f_m$ are real and simple whenever $m\geq m_0$, and
(2.2) shows that the same holds for $\phi(D)^mf$. \qed
\enddemo

\demo{Proof of Theorem 1.2} Again, we may assume that $f$ is monic
and $\phi(0)=1$. Let $d=\deg f$, and let $p, \alpha, \beta$ and the
polynomials $f_1, f_2, \dots$ be as in Theorem 2.1. We have
$\beta\ne 0$; and in the case where $p=2$ we must have $\beta>0$,
because we are assuming that $\phi''(0) - \phi'(0)^2\geq 0$. Hence
the corollary to Theorem 2.2 implies that $Z_C\left(\exp\left(\beta
D^p\right)M^d\right)>0$. By Theorem 2.1, $f_m \to \exp\left(\beta
D^p\right)M^d$ uniformly on compact sets in the complex plane. Hence
Rouche's theorem implies that there is a positive integer $m_0$ such
that $Z_C(f_m)>0$ whenever $m\geq m_0$, and (2.1) shows that the
same holds for $\phi(D)^mf$. \qed
\enddemo

In order to prove Theorem 2.1, we need some preliminaries. Let $\Bbb
C[x]$ denote the (complex) vector space of complex polynomials, let
$\Bbb C[x]^d$ denote the $(d+1)$-dimensional subspace of $\Bbb C[x]$
whose members are complex polynomials of degree $\leq d$, and let
$\|\phantom{f} \|_{\infty}$ denote the norm on  $\Bbb C[x]$ defined
by
$$
\|f\|_{\infty}= \max \{ |f^{(k)}(0)/k!| : 0\leq k\leq \deg f\}.
$$
Note that if $\langle f_m \rangle$ is a sequence of polynomials in
$\Bbb C[x]^d$, then $\|f_m\|_{\infty} \to 0$ if and only if $f_m \to
0$ uniformly on compact sets in the complex plane. When $\phi$ is a
formal power series (with complex coefficients) and $d$ is a
nonnegative integer, we denote the {\it operator norm} of
$\phi(D)|_{\Bbb C[x]^d}$ with respect to $\|\phantom{f} \|_{\infty}$
by $\|\phi(D) \|_{d}$, that is,
$$
\|\phi(D) \|_{d}=\sup\{\|\phi(D)f\|_{\infty}: f\in\Bbb C[x]^d \text{
and } \|f\|_{\infty}\leq 1 \}.
$$
If we denote the $d$-th partial sum of $\phi$ by $\phi|_d$, that is,
$$
\phi|_d(x)=\sum_{k=0}^d \frac{\phi^{(k)}(0)}{k!}x^k,
$$
then the restriction of $\phi(D)$ to $\Bbb C[x]^d$ is completely
determined by $\phi|_d$. Hence there are positive constants $A_d$
and $B_d$ such that
$$
A_d \|\phi(D) \|_{d}\leq \|\phi|_d\|_{\infty}\leq B_d\|\phi(D)
\|_{d}
$$
for all $\phi$.

For $c\ne 0$ we define the {\it dilation operator} $\Delta_c$ by
$$
\left(\Delta_c f\right)(x) = f(cx).
$$
It is then easy to see that
$$
\Delta_c\left(\phi(D)f\right)=\phi(c^{-1}D)(\Delta_c f)\qquad (c\ne
0),  \tag 2.3
$$
whenever $\phi$ is a formal power series and $f\in\text{\rm
dom\,}\phi(D)$.

\demo{Proof of Theorem 2.1} Let $r=\max\{p, d\}$. If $\tilde\phi$ is
a formal power series and $\tilde\phi|_r=\phi|_r$, then $\tilde\phi$
satisfies the identical assumptions in the theorem that are
satisfied by $\phi$, and we have $\tilde\phi(D)^mf=\phi(D)^mf$ for
all $m$. In other words, the theorem is about the first $r+1$
coefficients of $\phi$ only, and the coefficients
$\phi^{(n)}(0)/n!$, $n>r$, are irrelevant to the theorem. For this
reason, we may assume that $\phi^{(n)}(0)=0$ for all $n>r$. Then
there is a neighborhood $U$ of $0$ in the complex plane and there is
an analytic function $\psi$ in $U$ such that
$$
\log\phi(x)=\alpha x + \beta x^p + x^{p+1}\psi(x)\qquad (x\in U).
$$
We substitute $m^{-1/p}x$ for $x$ and multiply both sides by $m$ to
obtain
$$
m\log\phi\left(m^{-1/p}x\right)=m^{1-\frac{1}{p}}\alpha x + \beta
x^p +  m^{-1/p}x^{p+1}\psi\left(m^{-1/p}x\right)\qquad (x\in
m^{1/p}U).
$$
If we put
$$
\exp\left(-m^{1-\frac{1}{p}}\alpha
x\right)\phi\left(m^{-1/p}x\right)^m - \exp\left(\beta
x^p\right)=R_m(x),
$$
then $R_m$ is an entire function and we have
$$
R_m(x)=\exp\left(\beta
x^p\right)\left(\exp\left(m^{-1/p}x^{p+1}\psi\left(m^{-1/p}x\right)\right)-1\right)\qquad
(x\in m^{1/p}U).
$$
It is then clear that
$$
\sup_{|x|\leq R} |R_m(x)|=O\left(m^{-1/p}\right)\qquad (m\to\infty)
$$
for every $R>0$, and this implies that
$$
\left\|\exp\left(-m^{1-\frac{1}{p}}\alpha
D\right)\phi\left(m^{-1/p}D\right)^m - \exp\left(\beta
D^p\right)\right\|_d=O\left(m^{-1/p}\right) \tag 2.4
$$
as $m\to\infty$. Since $f$ is monic and of degree $d$, it follows
that
$$
\left\|m^{-d/p}\Delta_{m^{1/p}}f - M^d\right\|_{\infty}
=O\left(m^{-1/p}\right)\qquad (m\to\infty). \tag 2.5
$$
It is easy to see that (2.1) is equivalent to
$$
f_m=m^{-d/p}\exp\left(-m^{1-\frac{1}{p}}\alpha
D\right)\Delta_{m^{1/p}}\left(\phi(D)^mf\right);
$$
and (2.3) implies that the right hand side is equal to
$$
\exp\left(-m^{1-\frac{1}{p}}\alpha
D\right)\phi\left(m^{-1/p}D\right)^m\left(m^{-d/p}\Delta_{m^{1/p}}f
\right).
$$
Therefore we have
$$
\left\|f_m - \exp\left(\beta
D^p\right)M^d\right\|_{\infty}=O\left(m^{-1/p}\right)\qquad
(m\to\infty),
$$
by (2.4), (2.5) and the triangle inequality. This proves the
theorem. \qed
\enddemo

As we shall see soon, Theorem 2.2 is a consequence of known results
on Jensen polynomials and Mittag-Leffler functions. The following is
a simplified version of   \cite{6, Proposition 4.1}.

\proclaim{Proposition 2.3} Suppose that $\phi\in\Cal{LP}$,
 $q$ is a positive integer and   $f$ is given by
$$
f(x)=\sum_{k=0}^q \binom qk \phi^{(k)}(0) x^k.
$$
Suppose also that $\phi(0)\ne 0$ and $\phi$ is not of the form
$\phi(x)=p(x)e^{\alpha x}$, where $p$ is a polynomial and $\alpha\ne
0$. Then all the zeros of $f$ are real and simple.
\endproclaim

\remark{Remark} The polynomial $f$  is called the $q$-th {\it Jensen
polynomial} associated with $\phi$.
\endremark

The  {\it Mittag-Leffler functions} $E_1, E_2, \dots$ are given by
$$
E_p(x) = \sum_{k=0}^{\infty}\frac{x^k}{(pk)!}.
$$
It is  known that $E_p$ is of order $1/p$ and $E_p\in \Cal{LP}$ for
all $p$ (p.5 of \cite{14} and \cite {16, 19, 22}). Let $J_{(p,q)}$
denote the $q$-th Jensen polynomial associated with $E_p$:
$$
J_{(p,q)}(x)=\sum_{k=0}^q \binom qk E_p^{(k)}(0) x^k=\sum_{k=0}^q
\frac{q! x^k}{(q-k)!(pk)!}.
$$

\proclaim{Proposition 2.4} The zeros of $J_{(p,q)}$ are all negative
and simple for $p=2, 3, \dots$ and for $q=1, 2, \dots$.
\endproclaim

\demo{Proof} Suppose that $p\geq 2$ and $q\geq 1$. Then $E_p$ is of
order $\leq 1/2$, hence  it is not of the form $E_p(x)=p(x)e^{\alpha
x}$ where $p$ is a polynomial and $\alpha\ne 0$; and we have
$E_p(0)=1\ne 0$. Since $E_p\in\Cal{LP}$, Proposition 2.3 implies
that all the zeros of $J_{(p,q)}$ are real and simple. Finally, they
are all negative, because the coefficients of $J_{(p,q)}$ are all
positive. \qed\enddemo

\demo{Proof of Theorem 2.2} We have $d=pq+r$, $0\leq r\leq p-1$ and
$$
\left(\exp\left(-D^p\right)M^d\right)(x)=x^r \sum_{k=0}^q
\frac{(-1)^k d!}{k!(d-pk)!}x^{p(q-k)}.
$$
The right hand side is of the form $x^r f(x^p)$, where $f$ is a
monic polynomial of degree $q$ and $f(0)=(-1)^q d!/(q!r!)\ne 0$.
From this, we see that (1) is trivial, $\exp\left(-D^p\right)M^d$
has exactly $r$ zeros at the origin, and that the second assertion
of (2) follows from the first one. If $a\ne 0$ is a zero of
$\exp\left(-D^p\right)M^d$, then so are $e^{2k\pi i/p}a$, $k=0, 1,
\dots , p-1$, and they are distinct. Since
$\exp\left(-D^p\right)M^d$ has exactly $d=pq+r$ zeros in the whole
plane and has exactly $r$ zeros at the origin, it follows that
$\exp\left(-D^p\right)M^d$ has at most $q$ distinct positive zeros.
Hence it is enough to show that if $q\geq 1$, then
$\exp\left(-D^p\right)M^d$ has (at least) $q$ distinct positive
zeros.

Suppose that $q\geq 1$. We first consider the case where $d$ is a
multiple of $p$. In this case, we have $d=pq$, $r=0$ and
$$
\aligned \left(\exp\left(-D^p\right)M^d\right)(x) &=
 \sum_{k=0}^q
\frac{(-1)^k (pq)!}{k!(p(q-k))!}x^{p(q-k)}\\
&= \sum_{k=0}^q
\frac{(-1)^{q-k} (pq)!}{(q-k)!(pk)!}x^{pk}\\
&= (-1)^q\frac{(pq)!}{q!}\sum_{k=0}^q
\frac{q!}{(q-k)!(pk)!}(-x^p)^k\\
&= (-1)^q\frac{(pq)!}{q!}J_{(p,q)}(-x^p).
\endaligned
$$
Since $p\geq2$, Proposition 2.4 implies that all the zeros of
$J_{(p,q)}$ are negative and simple. Hence
$\exp\left(-D^p\right)M^d$ has exactly $q$ (=$\deg J_{(p,q)}$)
distinct positive zeros.

Finally, the result for the remaining case follows from an inductive
argument based on Rolle's theorem, because
$\left(\exp\left(-D^p\right)M^{pq+r}\right)(0)=0$ for $1\leq r\leq
p-1$,
$$
\exp\left(-D^p\right)M^d=\frac{1}{d+1}D\left(\exp\left(-D^p\right)M^{d+1}\right),
$$
and $\exp\left(-D^p\right)M^{p(q+1)}$ has exactly $q+1$ distinct
positive zeros.  \qed
\enddemo

\head 3. Proof of Theorem 1.3
\endhead

Let $\phi$ be a formal power series. First of all, we need to find a
sufficient condition for an entire function $f$ to be such that
$f\in\text{\rm dom\,} \phi(D)^m$  and  $\phi(D)^m f$ is not
identically equal to $0$ for all positive integers $m$. Let $\langle
C_n\rangle$ be a sequence of positive numbers. If
$|\phi^{(n)}(0)|<C_n$ for all $n$, we write $\phi \ll \langle
C_n\rangle$. More generally, if there are constants $c$ and $d$ such
that $c>0$, $d\geq 0$ and $\phi \ll \langle c(1+n)^d C_n\rangle$,
then we will write $\phi \prec \langle C_n\rangle$.

\proclaim{Lemma 3.1} Suppose that $\langle B_n\rangle$ is an
increasing sequence of positive numbers,
$$
B_mB_n\leq B_0B_{m+n} \qquad (m, n = 0, 1, 2, \dots),  \tag 3.1
$$
$\phi$ and $\psi$ are formal power series, $\phi, \psi \prec \langle
n!B_n\rangle$, $f$ is an entire function, and that  $f \prec \langle
(n!B_n)^{-1}\rangle$. Then $\phi \psi \prec \langle n!B_n\rangle$,
$f\in\text{\rm dom\,} \phi(D)$, $\phi(D) f \prec \langle
(n!B_n)^{-1}\rangle$ and $\psi(D)(\phi(D) f)= (\phi\psi)(D)f$.
\endproclaim

\demo{Proof} Suppose that $a, b$ are nonnegative constants, $\phi
\ll \langle (1+n)^a n! B_n\rangle$ and $\psi \ll \langle (1+n)^b n!
B_n\rangle$. Then
$$
\align \left|(\phi\psi)^{(n)}(0)\right|&\leq\sum_{k=0}^n\binom nk
\left|\phi^{(k)}(0)\right| \left|\psi^{(n-k)}(0)\right|\\
 &< B_0
(1+n)^{a+b+1} n! B_n \qquad (n=0, 1, 2, \dots),
\endalign
$$
hence $\phi \psi \prec \langle n!B_n\rangle$.

Now suppose that $c$ is a nonnegative constant, $f \ll \langle
(1+n)^c (n! B_n)^{-1}\rangle$, $R>0$, and $|x|\leq R$. Then
$$
\left| \frac{\phi^{(n)}(0)f^{(n+k)}(0)x^k}{n!k!}\right| \leq
\frac{B_0(1+n)^{a+c}(1+k)^c R^k}{n!(k!)^2 B_k},  \tag 3.2
$$
and we have
$$
\sum_{n, k\geq 0}\frac{B_0(1+n)^{a+c}(1+k)^c R^k}{n!(k!)^2 B_k} \leq
\sum_{n=0}^{\infty}\frac{(1+n)^{a+c}}{n!}\sum_{k=0}^{\infty}\frac{(1+k)^c
R^k}{(k!)^2 }<\infty.
$$
Hence the double series
$$
\sum_{n, k\geq 0}\frac{\phi^{(n)}(0)f^{(n+k)}(0)x^k}{n!k!}
$$
converges absolutely and uniformly on compact sets in the complex
plane. As a consequence, the series
$$
\sum_{n=0}^{\infty}\frac{\phi^{(n)}(0)}{n!}f^{(n)}
$$
converges uniformly on compact sets in the complex plane, that is,
$f\in\text{\rm dom\,} \phi(D)$. Furthermore, the absolute
convergence of the double series  implies that
$$
\phi(D)f(x)=\sum_{k=0}^{\infty}\sum_{n=0}^{\infty}
\frac{\phi^{(n)}(0)f^{(n+k)}(0)}{n!}\frac{x^k}{k!} \qquad (x\in\Bbb
C),
$$
from which we obtain
$$
(\phi(D)f)^{(k)}(0) = \sum_{n=0}^{\infty}
\frac{\phi^{(n)}(0)f^{(n+k)}(0)}{n!} \qquad (k=0, 1, 2, \dots),
$$
and the assumptions imply that
$$
\left|(\phi(D)f)^{(k)}(0)\right| < \frac{B_0 (1+k)^c}{k! B_k}
\sum_{n=0}^{\infty} \frac{(1+n)^{a+c}}{n!} \qquad (k=0, 1, 2,
\dots),
$$
hence we have $\phi(D) f \prec \langle (n!B_n)^{-1}\rangle$.

Finally, an estimate which is similar to (3.2) shows that the triple
series
$$
\sum_{m, n, k\geq
0}\frac{\psi^{(m)}(0)\phi^{(n)}(0)f^{(m+n+k)}(0)x^k}{m!n!k!}
$$
converges absolutely for every $x\in\Bbb C$, hence the last
assertion follows. \qed
\enddemo

\proclaim{Corollary} Suppose that $\phi, \psi, f$ and $\langle B_n
\rangle$ are as in Lemma 3.1, $\mu$ is a nonnegative integer,
$\phi(x)\psi(x)=x^{\mu}$, and that $f$ is transcendental. Then
$f\in\text{\rm dom\,} \phi(D)^m$  and  $\phi(D)^m f$ is not
identically equal to $0$ for all positive integers $m$.
\endproclaim

\demo{Proof} An inductive argument shows that $f\in\text{\rm dom\,}
\phi(D)^m$, $\phi(D)^m f\in\text{\rm dom\,} \psi(D)^m$, and that
$\psi(D)^m\left(\phi(D)^m f\right) = f^{(m\mu)}$ for all $m$. Since
$f$ is transcendental, $f^{(m\mu)}$ is not identically equal to $0$
for all $m$, hence the same is true for $\phi(D)^m f$. \qed
\enddemo

\proclaim{Lemma 3.2} Suppose that $\langle B_n\rangle$ and $\phi$
are as in Lemma 3.1, $f$ is an entire function, $\langle f_N\rangle$
is a sequence of entire functions, $f_N \ll \langle (n!
B_n)^{-1}\rangle$ for all $N$, and that $f_N\to f$ as $N\to\infty$
uniformly on compact sets in the complex plane. Then $f_1, f_2,
\dots, f \in \text{\rm dom\,} \phi(D)$ and  $\phi(D)f_N\to \phi(D)f$
as $N\to\infty$ uniformly on compact sets in the complex plane.
\endproclaim

\demo{Proof} First of all, Lemma 3.1 implies that $f_N\in \text{\rm
dom\,} \phi(D)$ for all $N$. Since $f_N\to f$  uniformly on compact
sets in the complex plane, and since
$$
|f_N^{(n)}(0)|<(n! B_n)^{-1}\qquad (N=1,2,\dots;\ n=0,1,2,\dots),
$$
it follows that
$$
|f^{(n)}(0)|\leq(n! B_n)^{-1}\qquad (n=0,1,2,\dots),
$$
hence $f\in \text{\rm dom\,} \phi(D)$, by Lemma 3.1.

To prove the uniform convergence on compact sets in the complex
plane, let $R>0$ and $\epsilon>0$ be arbitrary. Suppose that $a$ is
a nonnegative constant and $\phi\ll \langle (1+n)^a n! B_n\rangle$.
If we put
$$
b=\sum_{k=0}^{\infty}\frac{B_0 R^k}{(k!)^2 B_k},
$$
then it is easy to see that
$$
|f_N^{(n)}(x)|<\frac{b}{n!B_n} \qquad (|x|\leq R;\ N=1,2,\dots;\
n=0,1,2,\dots),
$$
and that
$$
|f^{(n)}(x)|\leq\frac{b}{n!B_n} \qquad (|x|\leq R;\ n=0,1,2,\dots).
$$
Let $\nu$ be a positive integer such that
$$
b\sum_{n=\nu + 1}^{\infty}\frac{(1+n)^a}{n!}<\epsilon.
$$
Then there is a positive integer $N_0$ such that
$$
\left|\sum_{n=0}^{\nu}\frac{\phi^{(n)}(0)}{n!}\left(f_N^{(n)}(x)-f^{(n)}(x)\right)\right|<\epsilon
\qquad (|x|\leq R;\ N\geq N_0),
$$
because $f_N\to f$  uniformly on compact sets in the complex plane.

Now, suppose that $|x|\leq R$ and $N\geq N_0$. Then we have
$$
\multline |\phi(D)f_N(x) - \phi(D)f(x)| \leq \\
\left|\sum_{n=0}^{\nu}\frac{\phi^{(n)}(0)}{n!}\left(f_N^{(n)}(x)-f^{(n)}(x)\right)\right|
+2\sum_{n=\nu+1}^{\infty}\frac{|\phi^{(n)}(0)|}{n!}\frac{b}{n!B_n}<3\epsilon.
\endmultline
$$
This completes the proof.  \qed

\enddemo

\proclaim{Corollary} Under the same assumptions as in Lemma 3.2,
$\phi(D)^m f_N \to \phi(D)^m f$ as $N\to\infty$ uniformly on compact
sets in the complex plane for every positive integer $m$.
\endproclaim

\demo{Proof} Lemma 3.1 implies that $\phi^m\prec\langle
n!B_n\rangle$ for all positive integers $m$.  \qed
\enddemo

We denote the open disk with center at $a$ and radius $r$ by
$D(a;r)$, and its closure by $\bar{D}(a;r)$. For a complex constant
$c$ we define the {\it translation operator} $T^c$ by $(T^c
f)(x)=f(x+c)$. It is clear that if $f$ is a monic polynomial of
degree $d$, then $c^{-d}T^c f \to 1$ as $|c|\to\infty$ uniformly on
compact sets in the complex plane. This observation leads to the
following:

\proclaim{Lemma 3.3} Suppose that $\phi$ is a formal power series,
$f$ and $g$ are polynomials, $a_1, \dots, a_N$ are zeros of
$\phi(D)f$, $b$ is a zero of $\phi(D)g$, and that neither $\phi(D)f$
nor $\phi(D)g$ is identically equal to $0$. Then for every $c\in\Bbb
C$ the polynomial  $\phi(D)(fT^c g)$ is not identically equal to
$0$, and for every $\epsilon>0$ there is an $R>0$ such that if
$|c|>R$, then $\phi(D)(fT^c g)$ has a zero in each of the disks
$D(a_1;\epsilon), \dots,D(a_N;\epsilon)$ and $D(b-c; \epsilon)$.
\endproclaim

\demo{Proof} The assumptions imply that neither $f$ nor $g$ is
identically equal to $0$. In particular, we have $\deg (fT^c g) \geq
\deg f$, hence the first assertion follows, because $\phi(D)f$ is
not identically equal to zero.

Let $\epsilon>0$.  We first observe that if $c$ is a constant, then
$\phi(D)(fT^c g)$ has a zero in $D(b-c; \epsilon)$ if and only if
$\phi(D)(gT^{-c} f)$ has a zero in $D(b; \epsilon)$. Since neither
$f$ nor $g$ is  identically equal to $0$, we may assume that $f$ and
$g$ are monic. Then $c^{-\deg g} fT^c g \to f$ and $(-c)^{-\deg f}
gT^{-c} f \to g$ as $|c|\to\infty$ uniformly on compact sets in the
complex plane. Hence there is an $R>0$ such that if $|c|>R$, then
$\phi(D)(fT^c g)$ has a zero in each of the disks $D(a_1;\epsilon),
\dots,D(a_N;\epsilon)$ and $\phi(D)(gT^{-c} f)$ has a zero in $D(b;
\epsilon)$.  \qed\enddemo

The following characterization of the class $\Cal{LP}$ given in
\cite{18, 21} will play a crucial role in the proof of Theorem 1.3.

\proclaim{Theorem (P\'olya)} Let $\phi$ be a formal power series
with real coefficients. Then $\phi\in\Cal{LP}$ if and only if
$Z_C(\phi(D)M^d)=0$ for all positive integers $d$.
\endproclaim

\proclaim{Corollary} Suppose that $\phi$ is a formal power series
with real coefficients and $\phi$ does not represent a
Laguerre-P\'olya function. Then there is a positive integer $d_0$
such that $Z_C(\phi(D)M^d)>0$ for all $d\geq d_0$.
\endproclaim

\demo{Proof} By  P\'olya's theorem, there is a  positive integer
$d_0$ such that $Z_C(\phi(D)M^{d_0})>0$, and Rolle's theorem implies
that if $Z_C(\phi(D)M^{d+1})=0$, then $Z_C(\phi(D)M^{d})=0$.  \qed
\enddemo

\demo{Proof of Theorem 1.3} We will construct a sequence $\langle
d(k)\rangle$ of positive integers and a sequence $\langle
\gamma(k)\rangle$ of positive numbers such that
$\sum_{k=1}^{\infty}d(k)\gamma(k)<\infty$ and the entire function
$f$ represented by
$$
f(x)=\prod_{k=1}^{\infty} \left(1+\gamma(k)x\right)^{d(k)}
$$
has the desired property.

Since $\phi$ does not represent a Laguerre-P\'olya function, it
follows that neither does the formal power series $\phi^m$   for
every positive integer $m$. Hence the corollary to P\'olya's theorem
implies that there is an increasing sequence $\langle d(m)\rangle$
of positive integers such that $Z_C(\phi(D)^m M^{d(m)})>0$ for all
positive integers $m$. Since $\langle d(m)\rangle$ is increasing, we
have $Z_C(\phi(D)^m M^{d(k)})>0$ whenever $m\leq k$. For each pair
$(m,k)$ of positive integers with $m\leq k$ choose a nonreal zero of
$\phi(D)^m M^{d(k)}$ in the upper half plane, denote it by $a(m,k)$
and set $r(m, k)=\I a(m,k)/2$. It is obvious that $r(m,k)>0$, and
that $\bar D(a(m,k)-\gamma; r(m,k))\cap\Bbb R=\emptyset$ for all
$\gamma\in\Bbb R$. The assumption also implies that
$\phi^{(n)}(0)\ne 0$ for some $n$, hence there is a nonnegative
integer $\mu$ and there is a formal power series $\psi$ such that
$\phi(x)\psi(x)=x^{\mu}$. Choose an increasing sequence $\langle A_n
\rangle$ of positive numbers such that $\phi, \psi \ll \langle A_n
\rangle$, and define $\langle B_n \rangle$ by $B_0=A_0$, $B_1=A_1$
and
$$
B_{n+1}=\max\left[\{A_{n+1}\} \cup \{ B_0^{-1}B_kB_{n+1-k}:
k=1,\dots, n\}\right] \qquad (n=1, 2, \dots).
$$
It is clear that $\langle B_n\rangle$ is an increasing sequence of
positive numbers, $\langle B_n\rangle$ satisfies (3.1), and that
$\phi, \psi \prec \langle n!B_n\rangle$.

For $k=1, 2, \dots$ and for $\gamma>0$ define $g_{k,\gamma}$ by
$$
g_{k,\gamma}(x)=\left(1+\gamma x \right)^{d(k)},
$$
that is, $g_{k,\gamma}=\gamma^{d(k)}T^{1/\gamma} M^{d(k)}$. From the
definition, it follows that $g_{k,\gamma}$ is a real polynomial of
degree $d(k)$,  $g_{k,\gamma}(0)=1$,  $\phi(D)^m g_{k,\gamma}$ is
not identically equal to $0$ for $1\leq m\leq k$,
$$
\left(\phi(D)^m g_{k,\gamma}\right)\left(a(m,k) -
\gamma^{-1}\right)=0 \qquad (1\leq m\leq k),  \tag 3.3
$$
and that $g_{k,\gamma} \to 1$ as $\gamma\to 0$ uniformly on compact
sets in the complex plane.

Since $g_{1,\gamma}(0)=1<2$ and $g_{1,\gamma}$ is a polynomial of
degree $d(1)$ for every $\gamma>0$, and since   $g_{1,\gamma}\to 1$
as $\gamma\to 0$ uniformly on compact sets in the complex plane,
there is a positive number $\gamma(1)$ such that
$g_{1,\gamma(1)}\ll\langle 2 B_0(n! B_n)^{-1}\rangle$. From the
definition, the polynomial $\phi(D)g_{1,\gamma(1)}$ is not
identically equal to $0$, and from (3.3) we have $\left(\phi(D)
g_{1,\gamma(1)}\right)\left(a(1,1) - \gamma(1)^{-1}\right)=0 $. Now
suppose that $\gamma(1), \dots, \gamma(N)$ are positive numbers,
$$
\prod_{k=1}^{N} g_{k, \gamma(k)}\ll\langle 2 B_0(n!
B_n)^{-1}\rangle, \tag 3.4
$$
and that for each $m\in\{1, \dots, N\}$ the closures of the disks
$$
D\left(a(m, k)-\gamma(k)^{-1} ; r(m,k)\right) \qquad (m\leq k\leq N)
\tag 3.5
$$
are mutually disjoint and the polynomial $\phi(D)^m\left(
\prod_{k=1}^{N} g_{k, \gamma(k)}\right)$ has a zero in each of these
disks. Suppose also that the polynomials $\phi(D)^m\left(
\prod_{k=1}^{N} g_{k, \gamma(k)}\right)$, $m=1, \dots, N$ are not
identically equal to $0$. Since $\prod_{k=1}^{N} g_{k, \gamma(k)}$
is a polynomial, $g_{N+1, \gamma}$ is a polynomial of degree
$d(N+1)$ for every $\gamma>0$, and since $g_{N+1, \gamma}\to 1$ as
$\gamma\to 0$ uniformly on compact sets in the complex plane, (3.4)
implies that there is a $\delta>0$ such that
$$
\left(\prod_{k=1}^{N} g_{k, \gamma(k)}\right)g_{N+1,
\gamma}\ll\langle 2 B_0(n! B_n)^{-1}\rangle \qquad(0<\gamma<\delta).
\tag 3.6
$$
From Lemma 3.3, it follows that for each $m\in\{1, \dots, N\}$ there
is an $R_m>0$ such that if $|c|>R_m$, then  $\phi(D)^m\left(\left(
\prod_{k=1}^{N} g_{k, \gamma(k)}\right)T^c M^{d(N+1)}\right)$ has a
zero in each of the disks given in (3.5) and also has a zero in the
disk $D\left(a(m, N+1) -c; r(m, N+1)\right)$, because
$\phi(D)^m\left( \prod_{k=1}^{N} g_{k, \gamma(k)}\right)$ has a zero
in each of the disks given in (3.5) and $\left(\phi(D)^m
M^{d(N+1)}\right)\left(a(m, N+1)\right)=0$. By taking $R_m$
sufficiently large, we may assume that
$$
\bar D\left(a(m, k)-\gamma(k)^{-1} ; r(m,k)\right)\cap \bar
D\left(a(m, N+1) -c; r(m, N+1)\right)=\emptyset
$$
for $|c|>R_m$ and for $m\leq k\leq N$. Since $\phi(D)^{N+1}
M^{d(N+1)}$ has a zero at $a(N+1, N+1)$ and $r(N+1, N+1)>0$, Lemma
3.3 implies that there is an $R_{N+1}>0$ such that if $|c|>R_{N+1}$,
then $\phi(D)^{N+1}\left(\left( \prod_{k=1}^{N} g_{k,
\gamma(k)}\right)T^c M^{d(N+1)}\right)$ has a zero in $D\left(a(N+1,
N+1) -c; r(N+1, N+1)\right)$. Let $\gamma(N+1)$ be such that
$0<\gamma(N+1)<\min\{ \delta, R_1^{-1}, \dots, R_{N}^{-1},
R_{N+1}^{-1}\}$. Then (3.6) implies that
$$
\prod_{k=1}^{N+1} g_{k, \gamma(k)}\ll\langle 2 B_0(n!
B_n)^{-1}\rangle .
$$
The construction shows that for each $m\in\{1, \dots, N+1\}$ the
closures of the disks
$$
D\left(a(m, k)-\gamma(k)^{-1} ; r(m,k)\right) \qquad (m\leq k\leq
N+1)
$$
are mutually disjoint and the polynomial $\phi(D)^m\left(
\prod_{k=1}^{N+1} g_{k, \gamma(k)}\right)$ has a zero in each of
these disks. Finally,  the polynomials $\phi(D)^m\left(
\prod_{k=1}^{N+1} g_{k, \gamma(k)}\right)$, $m=1, \dots,  N+1$, are
not identically equal to $0$, by Lemma 3.3.

By induction, this process produces a sequence $\langle
\gamma(k)\rangle$ of positive numbers which has the following
properties: \roster
\item For each positive integer $N$ we have
$$
\prod_{k=1}^{N} g_{k, \gamma(k)}\ll\langle 2 B_0(n!
B_n)^{-1}\rangle.
$$
\item For each positive integer $m$ the closed disks
$$
\bar D\left(a(m, k)-\gamma(k)^{-1} ; r(m,k)\right) \qquad (k=m, m+1,
m+2, \dots) \tag 3.7
$$
are mutually disjoint.
\item
For each positive integer $m$ the polynomial $\phi(D)^m\left(
\prod_{k=1}^{N} g_{k, \gamma(k)}\right)$ has a zero in each of the
disks given in (3.5), whenever $N\geq m$.
\endroster
For $N=1, 2, \dots$ we set $f_N=\prod_{k=1}^{N} g_{k, \gamma(k)}$,
 that is,
$$
f_N(x)=\prod_{k=1}^{N} \left(1 + \gamma(k)x \right)^{d(k)}.
$$
From (1), it follows that
$$
0\leq f_N^{(n)}(0)<2 B_0(n!B_n)^{-1} \qquad (N=1, 2, \dots ;\ n=0,
1, 2, \dots). \tag 3.8
$$
In particular, we have
$$
\sum_{k=1}^N d(k)\gamma(k)=f_N'(0)<2B_0/B_1 \qquad (N=1, 2, \dots),
$$
hence the infinite product $\prod_{k=1}^{\infty} \left(1 +
\gamma(k)x \right)^{d(k)}$ represents an entire function of genus
$0$. Let $f$ denote the entire function. It is then obvious that $f$
is transcendental, $f\in\Cal{LP}$, $f_N\to f$ uniformly on compact
sets in the complex plane, and that
$$
0<f^{(n)}(0)\leq 2 B_0(n!B_n)^{-1} \qquad ( n=0, 1, 2, \dots).
$$

To complete the proof, let $m$ be a positive integer. From the
corollary to Lemma 3.1, it follows that $f\in\text{\rm dom\,}
\phi(D)^m$  and  $\phi(D)^m f$ is not identically equal to $0$; and
from (3.8) and the corollary to Lemma 3.2, we see that $\phi(D)^m
f_N \to \phi(D)^m f$ as $N\to\infty$ uniformly on compact sets in
the complex plane. Furthermore, $\phi(D)^m f_N$ has a zero in each
of the disks given in (3.5) whenever $N\geq m$. Hence $\phi(D)^mf$
has a zero in each of the closed disks given in (3.7) which are
mutually disjoint and do not intersect the real axis. Therefore
$Z_C(\phi(D)^m f)=\infty$. \qed

\enddemo

\head 4. Some Consequences of Theorems 2.1 and 2.2
\endhead

In this brief section, we will be concerned with the asymptotic
behavior of the distribution of zeros of $\phi(D)^m f$ as
$m\to\infty$, in the case where the coefficients of $\phi$ are
complex numbers and $f$ is a complex polynomial.  When $f$ is an
entire function, we denote its zero set  by $\Cal Z(f)$, that is,
$\Cal Z(f)=\{z\in\Bbb C : f(z)=0\}$, and for $a\in\Cal Z(f)$ the
multiplicity by $\text{\rm{m}}(a,f)$.

Let $\phi$, $p$, $\alpha$, $\beta$, $f$, $d$ and $f_1, f_2, \dots$
be as in Theorem 2.1. Then $\beta\ne 0$ and $f_m\to \exp(\beta
D^p)M^d$ uniformly on compact sets in the complex plane. We also
have
$$
\Cal Z(\phi(D)^m f) = -m\alpha + m^{1/p}\Cal Z(f_m)  \tag 4.1
$$
and
$$
 \text{\rm m} (a, \phi(D)^m f)=\text{\rm m}
(m^{-1/p}(a+m\alpha), f_m)  \tag 4.2
$$
for all $a\in \Cal Z(\phi(D)^m f))$.  Let $\epsilon>0$ be so small
that the disks $D(b;\epsilon)$, $b\in\Cal Z\left(\exp(\beta
D^p)M^d\right)$, are mutually disjoint. Then Rouche's theorem
implies that there is a positive integer $m_0$ such that
$$
\sum_{c\in D(b;\epsilon)\cap \Cal Z(f_m)} \text{\rm{m}}(c,
f_m)=\text{\rm{m}}\left(b, \exp(\beta D^p)M^d\right)   \tag 4.3
$$
holds for all $b\in\Cal Z\left(\exp(\beta D^p)M^d\right)$ for all
and $m\geq m_0$. As a consequence, we have
$$
\Cal Z(f_m)\subset D(0;\epsilon) + \Cal Z\left(\exp(\beta
D^p)M^d\right) \tag 4.4
$$
for all $m\geq m_0$. Let $\gamma$ be a complex number such that
$\gamma^p = -\beta$. Then $\gamma\ne 0$ and  we have
$$
\Cal Z\left(\exp(\beta D^p)M^d\right)=\gamma \Cal Z\left(\exp(-
D^p)M^d\right),
$$
because
$$
\left(\exp(\beta D^p)M^d\right)(x)=\gamma^d \left(\exp(-
D^p)M^d\right)(x/\gamma).
$$
Now, (4.1) and (4.4) imply that
$$
\Cal Z(\phi(D)^m f) \subset -m\alpha + m^{1/p}\left(D(0;\epsilon) +
\gamma\Cal Z\left(\exp(- D^p)M^d\right)  \right)  \tag 4.5
$$
holds for all $m\geq m_0$.

With the aid of   Theorem 2.2, the above results give us some
information on the zeros of $\phi(D)^m f$ for large values of $m$.
From Theorem 2.2, it follows that
$$
\Cal Z\left(\exp(- D^p)M^d\right) \subset S_p,
$$
where
$$
S_p = \bigcup_{k=0}^{p-1}\left\{re^{2k\pi i/p} : r\geq 0\right\}.
$$
It  also  follows from Theorem 2.2 that if $d\equiv 0\text{ or
}1\mod p$, then all the zeros of $\exp(- D^p)M^d$ are  simple. Hence
(4.5) implies that for every $\epsilon>0$ there is a positive
integer $m_0$ such that
$$
\Cal Z(\phi(D)^m f) \subset -m\alpha + N(0, m^{1/p}\epsilon) +
\gamma S_p
$$
 for all $m\geq m_0$, and (4.1) through (4.3) imply that if
$d\equiv 0\text{ or }1\mod p$, then all the zeros of $\phi(D)^m f$
are  simple whenever $m$ becomes sufficiently large.

\enddocument